\documentclass[10pt,reqno]{article}
\usepackage{amssymb, amsmath, amsthm, amsfonts, amscd}
\usepackage[english]{babel}
\usepackage{graphicx, epsfig}

\newtheorem{theorem}{Theorem}[section]

\newcommand{\nm}{\noalign{\smallskip}}
\newcommand{\ds}{\displaystyle}

\newcommand{\p}{\partial}
\newcommand{\pd}[2]{\frac {\p #1}{\p #2}}
\newcommand{\eqnref}[1]{(\ref {#1})}

\newcommand{\Rbb}{\mathbb{R}}

\newcommand{\la}{\langle}
\newcommand{\ra}{\rangle}

\newcommand{\Kcal}{\mathcal{K}}

\newcommand{\Scal}{\mathcal{S}}

\newcommand{\Gvf}{\varphi}

\newcommand{\Gl}{\lambda}

\newcommand{\Gs}{\sigma}

\newcommand{\Gz}{\zeta}
\newcommand{\GD}{\Delta}

\newcommand{\GO}{\Omega}

\newcommand{\beq}{\begin{equation}}
\newcommand{\eeq}{\end{equation}}

\numberwithin{equation}{section}
\numberwithin{figure}{section}

\begin{document}
\title{Coated inclusions of finite conductivity neutral to multiple fields in two dimensional conductivity or anti-plane elasticity\thanks{\footnotesize This work is supported by Korean Ministry of Education, Sciences and Technology through NRF grants Nos. 2010-0004091 and 2010-0017532.}}

\author{Hyeonbae Kang\thanks{Department of Mathematics, Inha University, Incheon
402-751, Korea (hbkang@inha.ac.kr, hdlee@inha.ac.kr).}  \and Hyundae Lee\footnotemark[2]}

\date{}
\maketitle

\begin{abstract}
We consider the problem of neutral inclusions for two-dimensional conductivity and anti-plane elasticity. The neutral inclusion, when inserted in a matrix having a uniform field, does not disturb the field outside the
inclusion. The inclusion consists of a core and a shell. We show that if the inclusion is neutral to two linearly independent fields, then the core and the shell are confocal ellipses.
\end{abstract}

\section{Introduction}

The purpose of this paper is to prove that if a coated inclusion in two dimensions is neutral to multiple uniform fields, then the core and the shell are confocal ellipses.

We begin by defining the neutral inclusions. Let $D$ and $\GO$ be bounded domains with Lipschitz boundaries in $\Rbb^2$ such that $\overline{D} \subset \GO$ so that $D$ is the core and $\GO \setminus D$ is the shell. The conductivity is $\Gs_c$ in the core, $\Gs_s$ in the shell, and $\Gs_m$ in the matrix ($\Rbb^2 \setminus \GO$). So the conductivity distribution is given by
$$
\Gs = \Gs_c \chi(D) + \Gs_s \chi(\GO \setminus D) + \Gs_m\chi(\Rbb^2 \setminus \GO)
$$
where $\chi$ is the characteristic function. The conductivities $\Gs_c$ and $\Gs_s$ are assumed to be isotropic (scalar), but we allow $\Gs_m$ to be anisotropic, {\it i.e.}, a positive definite symmetric constant matrix. For a given function $h$ with $\nabla\cdot\Gs_m \nabla h=0$ in $\Rbb^2$, we consider
\beq\label{trans}
\left\{
\begin{array}{ll}
\nabla \cdot \Gs \nabla u = 0 \quad &\mbox{in } \Rbb^2, \\
u(x)-h(x) = O(|x|^{-1}) \quad &\mbox{as } |x| \to \infty,
\end{array}
\right.
\eeq
This problem can be regarded as a conductivity problem or an anti-plane elasticity problem.

The inclusion $\GO$ (or $\Gs$) is said to be neutral to the field $-\nabla h$ if the solution $u$ to \eqnref{trans} satisfies $u(x)-h(x)=0$ in $\Rbb^2 \setminus \GO$. So the neutral inclusion to the field $-\nabla h$ does not perturb the field outside the inclusion. In the imaging point of view, it means that the field $-\nabla h$ can not probe the neutral inclusion.

A particular interest lies in the inclusions neutral to uniform fields, {\it i.e.}, $h(x)=a \cdot x$ for some constant vector $a$. If $D$ and $\GO$ are concentric disks, say $D=\{\, |x|< r_1\,\}$ and $\GO=\{\, |x|< r_2\,\}$, and $\Gs_m$ is isotropic, then one can easily see that $\GO$ is neutral to all uniform fields if the following relation holds:
\beq\label{neutral_cond}
(\Gs_s + \Gs_c)(\Gs_m - \Gs_s) + f (\Gs_s - \Gs_c)(\Gs_m + \Gs_s) =0
\eeq
where $f=r_1^2/r_2^2$ (the volume fraction). Much interest in neutral inclusions (to uniform fields) was aroused by the work of Hashin \cite{hashine}. He showed that since insertion of neutral inclusions does not perturb the outside uniform field, the effective conductivity of the assemblage filled with coated inclusions of many different scales is $\Gs_m$ satisfying \eqnref{neutral_cond}. It is also proved that this effective conductivity is a bound of the Hashin-Shtrikman bounds of the effective conductivity of arbitrary two phase composite. We refer to a book of Milton \cite{milton} for development on neutral inclusions in relation to theory of composites.

Another interest in neutral inclusions has aroused in relation to imaging and invisibility cloaking by transformation optics. In this regard, we first observe that in general the solution $u$ to \eqnref{trans} satisfies $u(x)-h(x)=O(|x|^{-1})$ as $|x| \to \infty$. But, if the inclusion is neutral to all uniform fields, then the linear part of $h$ is unperturbed and one can show using multi-polar expansions that $u(x)-h(x)=O(|x|^{-2})$ as $|x| \to \infty$ for any $h$ (not necessarily linear). It means that it is harder to probe the neutral inclusions using $u-h$.  Recently, Ammari {\it et al} \cite{AKLL1} extend the idea of neutral inclusions to construct multi-coated circular structures which are neutral not only to uniform fields but also to fields of higher order up to $N$ for a given integer $N$ (which is called the GPT-vanishing structure, GPT for generalized polarization tensor), so that the solution $u$ to \eqnref{trans} satisfies $u(x)-h(x)=O(|x|^{-N-1})$ as $|x| \to \infty$ for any $h$. This structure has a strong connection to the cloaking by transformation optics. The transformation optics proposed by Pendry {\it et al} transforms a punctured disk (sphere) to an annulus to achieve perfect cloaking (the same transform was used to show non-uniqueness of the Calder\'on's problem by Greenleaf {\it et al} \cite{glu}). Kohn {\it et al} \cite{kohn1} showed that if one transforms a disk with a small hole, then one can avoid singularities of the conductivity which occur on the inner boundary of the annulus and achieve near-cloaking instead of perfect cloaking. In \cite{AKLL1} it is shown that if we coat the hole by the GPT-vanishing structure before transformation, then near-cloaking is dramatically enhanced. See also \cite{AGJKLL, AKLL2, AKLLY, kohn2} for further development to Helmholtz and Maxwell's equations.

There have been some work on neutral inclusions of general shape. Non-elliptic inclusions neutral to a single uniform field have been constructed by Milton and Serkov \cite{MS} when the conductivity $\Gs_c$ of the core is either $0$ or $\infty$, and by Jarczyk and Mityushev \cite{JM} when $\Gs_c$ is finite. In \cite{MS}, it is also proved that if an inclusion is neutral to all uniform field (or equivalently, to two linearly independent uniform fields), then the core and shell are confocal ellipses, when $\Gs_c$ is $0$ or $\infty$. That the confocal ellipses (and ellipsoids) are neutral to all uniform fields was proved by Kerker \cite{kerker}. 

The purpose of this paper is to prove that confocal ellipses are the only inclusions which are neutral to two fields even when $\Gs_c$ is finite. We emphasize that the method of \cite{MS} can not be applied to the case when $\Gs_c$ is finite. There the conformal mapping from the shell ($\GO \setminus D$) onto an annulus was used. But, there is no conformal mapping from $\GO$ onto a disk which maps $D$ to a concentric disk, except for some very special cases. 

More precisely we prove the following theorem. We emphasize that the theorem holds for any conductivity $\Gs_c$ of the core (finite, $0$, or $\infty$).
\begin{theorem}\label{mainthm}
Let $D$ and $\GO$ be bounded domains with Lipschitz boundaries in $\Rbb^2$ such that $\overline{D} \subset \GO$. If $\GO$ is neutral to $-\nabla x_j$ for $j=1,2$, then $D$ and $\GO$ are confocal ellipses.
\end{theorem}

The key observation in proving Theorem \ref{mainthm} is that if the inclusion is neutral to two fields, then the fields inside the core is also uniform. Using this fact, we are able to set up a free boundary value problem. We then use a conformal mapping to show that the solution to the free boundary value problem is a pair of confocal ellipses.

This paper is organized as follows: In the next section we review layer potential representation of the solution to \eqnref{trans}. In section 3, we show that if the inclusion is neutral to two fields, then the fields inside the core is uniform, and derive the free boundary value problem. In the last section we show that the solution to the free boundary value problem is confocal ellipses, and hence prove Theorem \ref{mainthm}.

\section{Layer potential representations of solutions}

Let $B$ be a bounded domain in $\Rbb^2$ with the Lipschitz boundary. The single layer potential $\Scal_{\p B} [\Gvf]$ of a density function $\Gvf
\in L^2(\p B)$ is defined by
\beq
\Scal_{\p B} [\Gvf] (x) := \frac{1}{2\pi} \int_{\p B} \ln |x-y| \Gvf (y) \, d\sigma(y) \;, \quad x \in \Rbb^2 .
\eeq
It satisfies the jump relation
\beq\label{singlejump}
\frac{\p}{\p\nu} \Scal_{\p B} [\Gvf] \big |_\pm (x) = \biggl( \pm \frac{1}{2} I + \Kcal_{\p B}^* \biggr) [\Gvf] (x),
\quad x \in \p B\;,
\eeq
where the operator $\Kcal_{\p B}$ is defined by
\beq\label{introkd}
\Kcal_{\p B} [\Gvf] (x) =
\frac{1}{2\pi} \int_{\p  B} \frac{\la
 y -x, \nu_y \ra}{|x-y|^2} \Gvf(y)\,d\sigma(y)\;, \quad x \in \p  B,
\eeq
and $\Kcal_{\p B}^*$ is its $L^2$-adjoint. The subscripts $+$ and $-$ indicate the limits from outside and inside $B$, respectively.

Suppose that $\Gs_m$ is isotropic. It is known (see, for example, \cite{book2}) that there is a pair $(\Gvf,\psi)\in L_0^2(\p D)\times L_0^2(\p \Omega)$ (the subscript $0$ indicates the mean value zero) such that the solution $u$ to \eqnref{trans} can be represented as
\begin{equation}\label{solrep}
u(x)= h(x) + \Scal_{\p D} [\Gvf](x) + \Scal_{\p \GO}[\psi](x), \quad x \in \Rbb^2.
\end{equation}
The transmission conditions (continuity of the potential and the flux) on the interfaces $\p D$ and $\p\GO$ are equivalent to
\begin{align}
(\Gl I - \Kcal_{\p D}^* ) [\Gvf] - \pd{}{\nu_{\p D}} \Scal_{\p \GO}[\psi] &= \pd{h}{\nu_{\p D}} \quad \mbox{on } \p D,\label{int_eq2}\\
- \pd{}{\nu_{\p \GO}} \Scal_{\p D}[\Gvf] + (\mu I - \Kcal_{\p \GO}^* ) [\psi] &= \pd{h}{\nu_{\p \GO}} \quad \mbox{on } \p \GO, \label{int_eq1}
\end{align}
where
\beq
\Gl= \frac{\sigma_c + \sigma_s}{2(\sigma_c - \sigma_s)}, \quad \mu=\frac{\sigma_s+ \sigma_m}{2(\sigma_s-\sigma_m)}.
\eeq
Here and throughout this paper $\pd{}{\nu_{\p D}}$ denotes the outward normal derivative on $\p D$. The system of integral equations \eqnref{int_eq2} and \eqnref{int_eq1} has a unique solution (see \cite{book2}). We emphasize that the representation \eqnref{solrep} is valid even if $\Gs_c$ is $0$ or $\infty$, and the solution $u$ is defined in $D$. If $\Gs_c=0$ then $\Gl=-1/2$, and if $\Gs_c=\infty$ then $\Gl=1/2$.

Note from the jump relation \eqnref{singlejump} that
\begin{align*}
\Gvf &= \pd{}{\nu_{\p D}} \Scal_{\p D}[\Gvf] \Big|_{+} - \pd{}{\nu_{\p D}} \Scal_{\p D}[\Gvf] \Big|_{-} \\
& =\pd{}{\nu_{\p D}} ( h + \Scal_{\p D} [\Gvf] + \Scal_{\p \GO}[\psi] ) \Big|_{+}
   - \pd{}{\nu_{\p D}}( h + \Scal_{\p D} [\Gvf] + \Scal_{\p \GO}[\psi] ) \Big|_{-} \\
&= \pd{u}{\nu_{\p D}} \Big|_{+} - \pd{u}{\nu_{\p D}} \Big|_{-} \\
&= \left( \frac{\Gs_c}{\Gs_s} -1 \right) \pd{u}{\nu_{\p D}} \Big|_{-}.
\end{align*}
So, we have
\beq\label{innersol}
\Gvf = \frac{2}{2\Gl-1} \pd{u}{\nu_{\p D}} \Big|_{-}.
\eeq
Similarly one can show that
\beq\label{psiform}
\psi= \left( 1- \frac{\Gs_m}{\Gs_s} \right) \pd{u}{\nu_{\p \GO}} \Big|_{+} = \frac{2}{2\mu+1} \pd{u}{\nu_{\p \GO}} \Big|_{+} .
\eeq
So, if $\sigma$ is neutral to the field $-\nabla h$, then $u=h$ in $\mathbb{R}^2\setminus \overline{\GO}$, and hence we have
\begin{equation}\label{psi_eq}
\psi = \frac{2}{2\mu+1} \pd{h}{\nu_{\p\GO}}.
\end{equation}

\section{Uniformity of the field in the core}

 Suppose that $\Gs$ is neutral to both $-\nabla x_1$ and $-\nabla x_2$. After a rotation if necessary, we may assume that $\Gs_m$ is given by
\beq
\Gs_m=\begin{bmatrix} \Gs_m^1 & 0 \\ 0 & \Gs_m^2 \end{bmatrix}.
\eeq
Let $u_j$, $j=1,2$, be the solution to \eqnref{trans} with $h(x)=x_j$. Then $u_j$ is the solution to \eqnref{trans} with $\Gs_m$ replaced by $\Gs_m^j$. So, we may represent it as
\beq\label{repreuj}
u_j(x)=x_j + \Scal_{\p D} [\Gvf_j](x) + \Scal_{\p \GO}[\psi_j](x).
\eeq
Let $n=(n_1, n_2)$ denote the outward unit normal to $\p \GO$ or $\p D$ and let
\beq
\mu_j=\frac{\sigma_s+ \sigma_m^j}{2(\sigma_s-\sigma_m^j)}, \quad j=1,2.
\eeq
Then we have from \eqnref{int_eq2} and \eqnref{int_eq1}
\begin{align}
(\Gl I - \Kcal_{\p D}^* ) [\Gvf_j] - \pd{}{\nu_{\p D}} \Scal_{\p \GO}[\psi_j] &= n_j \quad \mbox{on } \p D,\label{int_eq3}\\
- \pd{}{\nu_{\p \GO}} \Scal_{\p D}[\Gvf_j] + (\mu_j I - \Kcal_{\p \GO}^* ) [\psi_j] &= n_j \quad \mbox{on } \p \GO, \label{int_eq4}
\end{align}
and from \eqnref{psi_eq}
\begin{equation}\label{psi_eq2}
\psi_j = \frac{2}{2\mu_j+1} n_j \quad \mbox{on } \p\GO, \ j=1,2.
\end{equation}

Let $u_1^\perp$ be the harmonic conjugate of $u_1$ in each connected component of $\mathbb{R}^2\setminus (\p D \cup \p \GO)$. We may choose $u_1^\perp =x_2$ in $\Rbb^2 \setminus \GO$. We emphasize that even though $\GO \setminus D$ is not simply connected, a harmonic conjugate of $u_1$ exists there. In fact, a harmonic conjugate of $x_1 + \Scal_{\p \GO}[\psi_1]$ exists in $\GO$, and since $\Gvf_1 \in L^2_0(\p D)$ and $\Scal_{\p D} [\Gvf_1](x) = O(|x|^{-1})$, a harmonic conjugate of $\Scal_{\p D} [\Gvf_1]$ exists in $\Rbb^2 \setminus \overline{D}$.  Define $v$ by
\begin{equation}
v=\begin{cases}
\sigma_m^1 u_1^\perp \quad &\mbox{in}~\mathbb{R}^2\setminus \overline{\GO},\\
\sigma_s u_1^\perp +c_1\quad &\mbox{in}~ \GO\setminus \overline{D},\\
\sigma_c u_1^\perp +c_2 \quad &\mbox{in}~ D.\\
\end{cases}
\end{equation}
Here constants $c_1$ and $c_2$ are chosen so that $v$ is continuous across $\p \GO$ and $\p D$. Then $v$ satisfies
\begin{equation}\label{modified_eq}
\begin{cases}
\nabla \cdot (\widetilde{\sigma}\nabla v) = 0 \quad &\mbox{in}~\mathbb{R}^2,\\
v(x) - x_2 =0 &\mbox{in } \Rbb^2 \setminus \overline{\GO},
\end{cases}
\end{equation}
where \begin{equation} \widetilde{\sigma}=
\begin{cases}
(\sigma_m^1)^{-1}  \quad &\mbox{in}~\mathbb{R}^2\setminus \overline{\GO},\\
\sigma_s^{-1} \quad &\mbox{in}~ \GO\setminus \overline{D},\\
\sigma_c^{-1}  \quad &\mbox{in}~ D.\\
\end{cases}
\end{equation}
So $\widetilde{\sigma}$ is neutral to $-\nabla x_2$.

We represent $v$ as
$$
v(x) =x_2 + \Scal_{\p D} [\Gvf_3](x) + \Scal_{\p \GO}[\psi_3](x).
$$
Since
$$
\frac{\Gs_c^{-1} + \Gs_s^{-1}}{2(\Gs_c^{-1} - \Gs_s^{-1})} = -\Gl, \quad \frac{\Gs_s^{-1} +(\sigma_m^1)^{-1}}{2(\Gs_s^{-1}-(\sigma_m^1)^{-1})} = -\mu_1,
$$
the pair of potential $(\Gvf_3, \psi_3)$ satisfies
\begin{align}
(-\Gl I - \Kcal_{\p D}^* ) [\Gvf_3] - \pd{}{\nu_{\p D}} \Scal_{\p \GO}[\psi_3] = n_2 \quad & \mbox{on } \p D, \label{int_eq6}\\
 - \pd{}{\nu_{\p \GO}} \Scal_{\p D}[\Gvf_3] + (-\mu_1 I - \Kcal_{\p \GO}^* ) [\psi_3] = n_2  \quad & \mbox{on } \p \GO, \label{int_eq7}
\end{align}
and
\beq
\psi_3 = \frac{2}{-2\mu_1+1} n_2. \label{psi_eq3}
\eeq

It follows from \eqnref{psi_eq2} and \eqnref{psi_eq3} that
$$
\psi_2 = \frac{-2\mu_1 + 1}{2\mu_2 + 1} \psi_3,
$$
and since $\Scal_{\p D} [\Gvf_j] + \Scal_{\p \GO}[\psi_j]=0$ outside $\GO$ for $j=2,3$, we have
\begin{equation}
\Scal_{\p D}[\Gvf_2](x)= -\Scal_{\p \GO}[\psi_2](x)= -\frac{-2\mu_1 + 1}{2\mu_2 + 1} \Scal_{\p \GO}[\psi_3](x)=\frac{-2\mu_1 + 1}{2\mu_2 + 1} \Scal_{\p D}[\Gvf_3](x), \quad x \in \Rbb^2 \setminus \GO .
\end{equation}
Since $\Scal_{\p D}[\Gvf_2]$ and $\Scal_{\p D}[\Gvf_3]$ are harmonic in $\Rbb^2 \setminus \overline{D}$, we have
$$
\Scal_{\p D}[\Gvf_2](x) = \frac{-2\mu_1 + 1}{2\mu_2 + 1} \Scal_{\p D}[\Gvf_3](x), \quad x \in \Rbb^2 \setminus \overline{D}.
$$
Note that if $\Scal_{\p D}[\Gvf]=0$ outside $D$, then $\Gvf \equiv 0$. So we conclude
\begin{equation}\label{phi1_phi2}
\Gvf_2 =  \frac{-2\mu_1 + 1}{2\mu_2 + 1} \Gvf_3.
\end{equation}

We now see that \eqnref{int_eq6} can be written as
$$
(-\Gl I - \Kcal_{\p D}^* ) [\Gvf_2] - \pd{}{\nu_{\p D}} \Scal_{\p \GO}[\psi_2] = \frac{-2\mu_1 + 1}{2\mu_2 + 1} n_2 \quad  \mbox{on } \p D.
$$
By comparing this formula with \eqnref{int_eq3}, we deduce
\begin{equation}\label{gvf2}
\Gvf_2= \frac{\mu_1+\mu_2}{\Gl (2\mu_2+1)} n_2.
\end{equation}
We then have from \eqnref{innersol}
$$
\pd{u_2}{\nu} \Big|_{-} = \frac{2\Gl-1}{2} \Gvf_2= \frac{ (2\Gl-1)(\mu_1+\mu_2)}{2\Gl (2\mu_2+1)} n_2,
$$
and hence
\beq\label{u2x}
u_2(x)= \frac{ (2\Gl-1)(\mu_1+\mu_2)}{2\Gl (2\mu_2+1)} x_2 + c, \quad x \in D,
\eeq
for some constant $C$. So $-\nabla u_2$ is uniform in $D$.

By substituting \eqnref{psi_eq2} and \eqnref{gvf2} into \eqnref{repreuj}, we have
\beq
u_2(x) =x_2 + \Scal_{\p D} \left[ \frac{\mu_1+\mu_2}{\Gl (2\mu_2+1)} n_2 \right](x) + \Scal_{\p \GO} \left[ \frac{2}{2\mu_2+1} n_2 \right](x).
\eeq
It then follows from \eqnref{u2x} that
\beq\label{formn2}
(\mu_1+\mu_2) \Scal_{\p D} [n_2](x) + 2\Gl \Scal_{\p \GO} [n_2](x) =
\left\{
\begin{array}{ll}
0, \quad & x \in \Rbb^2 \setminus \GO, \\
\nm
\ds \left [\Gl(\mu_1-\mu_2-1)- \frac{\mu_1+\mu_2}{2} \right] x_2 +c_2, \quad & x \in D.
\end{array}
\right.
\eeq

By the exactly same argument with switched roles of $u_1$ and $u_2$, one can show that
\beq\label{formn1}
(\mu_1+\mu_2)  \Scal_{\p D} [n_1](x) + 2\Gl \Scal_{\p \GO} [n_1](x) =
\left\{
\begin{array}{ll}
0, \quad & x \in \Rbb^2 \setminus \GO, \\
\nm
\ds  \left[\Gl(\mu_2-\mu_1-1)- \frac{\mu_1+\mu_2}{2} \right]  x_1 + c_1, \quad & x \in D.
\end{array}
\right.
\eeq

Let $N_D(x)$ be the Newtonian potential on $D$, {\it i.e.},
\beq
N_D(x) := \frac{1}{2\pi} \int_{D} \ln |x-y| \, dy, \quad x \in \Rbb^2 \; .
\eeq
Then, we have
\beq
\Scal_{\p D} [n_j](x) = -\frac{\p}{\p x_j} N_D(x), \quad x \in \Rbb^2.
\eeq
It then follows from \eqnref{formn2} and \eqnref{formn1} that
$$
(\mu_1+\mu_2) N_D(x) + 2\Gl N_\GO(x) =
\left\{
\begin{array}{ll}
C_1, \quad & x \in \Rbb^2 \setminus \GO, \\
\nm
\ds   d_1 x_1^2+d_2x_2^2 - (c_1 x_1 + c_2 x_2) + C_2, \quad & x \in D,
\end{array}
\right.
$$
where $C_1$ and $C_2$ are constants and
\begin{align*} d_1=\frac{2\Gl(\mu_1-\mu_2+1)+(\mu_1+\mu_2)}{4},\\ d_2=\frac{2\Gl(\mu_2-\mu_1+1)+(\mu_1+\mu_2)}{4}.\end{align*} We may assume $c_1=c_2=0$ by translating $D$ and $\GO$. Observe that
$(\mu_1+\mu_2) N_D(x) + 2\Gl N_\GO(x)$ behaves as $\frac{(\mu_1+\mu_2) |D| + 2\Gl |\GO|}{2\pi} \ln |x|$ as $|x| \to \infty$, where $|D|$ denotes the area of $D$. So, we have
$$
(\mu_1+\mu_2) |D| + 2\Gl |\GO|=0,
$$
or
\beq
\frac{2\Gl}{\mu_1+\mu_2} = -f
\eeq
where $f=|D|/|\GO|$ (the volume fraction). It is worth mentioning that this relation shows that $\Gl$ and $\mu_1+\mu_2$ have opposite signs.
We finally have
\beq\label{newton}
N_D(x) -f N_\GO(x) =
\left\{
\begin{array}{ll}
C_1, \quad & x \in \Rbb^2 \setminus \GO, \\
\nm
\ds \frac{1-f(1+\mu_1-\mu_2)}{4} x_1^2 +\frac{1-f(1+\mu_2-\mu_1)}{4} x_2^2 + C_2, \quad & x \in D.
\end{array}
\right.
\eeq

It is worth mentioning that if the Newtonian potential of a simply connected domain is quadratic inside the domain, then the domain is an ellipse. This fact was proved by Dive \cite{div31} and Nikliborc \cite{nik32}, and is the key ingredient in resolution of the Eshelby's conjecture and the P\'olya-Szeg\"o conjecture by Kang and Milton \cite{KM} (see also \cite{kang}). In the next section we show that \eqnref{newton} implies that $D$ and $\GO$ are confocal ellipses.

\section{The free boundary value problem}

Let
$$
w(x) := \frac{f}{2} |x|^2 + 2(N_D(x) -f N_\GO(x)), \quad x \in \GO \setminus D.
$$
Note that $\GD N_\GO(x)=1$ for $x \in \GO$, $\GD N_\GO(x)=0$ for $x \notin \GO$, and $\nabla N_\GO$ is continuous across $\p \GO$. So, by \eqnref{newton} $w$ satisfies
\beq\label{free}
\left\{
\begin{array}{ll}
\ds \GD w= 0 \quad &\mbox{in } \GO \setminus \overline{D}, \\
\nm
\ds \nabla w = f x \quad &\mbox{on } \p\GO, \\
\nm
\ds \nabla w = x+(\mu_1-\mu_2)(x_1, -x_2) \quad &\mbox{on } \p D.
\end{array}
\right.
\eeq

We now show that if the problem \eqnref{free} admits a solution for some $f$, then $D$ and $\GO$ are confocal ellipses. Let $z=x_1+ix_2$ and let $w$ be a solution to \eqnref{free}. Define
$$
g(z): = 2 \p w(z) =  \frac{\p w}{\p x_1} - i\frac{\p w}{\p x_2}  .
$$
Then $g$ is holomorphic in $\GO \setminus \overline{D}$, $g(z)= f\overline{z}$ on $\p\GO$, and $g(z)=(\mu_1-\mu_2)z+ \overline{z}$ on $\p D$. Let $\Phi(\Gz)$ be the conformal mapping from $A=\{\Gz : 1<|\Gz|<r_0\}$ for some $r_0$ onto $\Omega\setminus \overline{D}$. Since $\p\GO$ and $\p D$ are Lipschitz, $\Phi$ is continuous on $\overline{A}$.
 Let $h(\Gz):=g(\Phi(\Gz))$. Then $h$ is holomorphic in $A$ and satisfies $h(\Gz)= f\overline{\Phi(\Gz)}$ on $|\Gz|=r_0$ and $h(\Gz)= \overline{\Phi(\Gz)}+(\mu_1-\mu_2)\Phi(\Gz)$ on $|\Gz|=1$.

Suppose that $\Phi$ admits the Laurent series expansion
\beq  \Phi(\Gz)= \sum_{n=-\infty}^{\infty} a_n \Gz^n. \eeq
Then we have
\beq
h(\Gz)=\begin{cases}
\ds f\sum_{n=-\infty}^{\infty} \overline{a_{-n}}r_0^{-2n} \Gz^{n}\quad &\mbox{on}~|\Gz|=r_0,\\
\ds \sum_{n=-\infty}^{\infty} \overline{a_{-n}}\Gz^{n}+(\mu_1-\mu_2)\sum_{n=-\infty}^{\infty} a_n \Gz^n\quad &\mbox{on}~|\Gz|=1.
\end{cases}
\eeq
Since the Laurent series expansion of $h$ is unique, we have
\beq \label{eq45}
f\overline{a_{-n}}r_0^{-2n}=\overline{a_{-n}}+(\mu_1-\mu_2)a_n
\eeq
for any integer $n$. Replacing $n$ with $-n$ and taking complex conjugates we also have
\beq   fa_{n}r_0^{2n}=a_{n}+(\mu_1-\mu_2)\overline{a_{-n}}. \eeq
Above two identities imply
\beq  [(1-fr_0^{-2n})(1-fr_0^{2n})-(\mu_1-\mu_2)^2]a_n=0 \eeq
for any integer $n$.
Since $(1-fr_0^{-2n})(1-fr_0^{2n})$ takes different values for different positive integers $n$, we know that there are only one positive $n$ with $a_n \ne 0$. By \eqnref{eq45}, the possibly nonzero coefficients of $\Phi(w)$ are $a_n$ and $a_{-n}$ and the univalence of $\Phi$ implies $n=1$. Finally, $\p D$ and $\p\GO$ are confocal ellipses since they are images of concentric circles by the map $\Phi(\Gz)=a_1 \Gz +a_{-1} \Gz^{-1}$.

\section*{Acknowledgement}

We are grateful to Habib Ammari and Graeme Milton for very helpful discussions on neutral inclusions. We also thank G. Milton for pointing out the reference \cite{kerker} to us.

\end{document}